\documentclass[12pt,draft]{amsart}
\usepackage{amsmath,amsthm,latexsym,amscd,amsbsy,amssymb,fleqn,leqno}
\chardef\bslash=`\\ 

\makeatletter
\def\verbatim{\interlinepenalty\@M \@verbatim
  \leftskip\@totalleftmargin\advance\leftskip2pc
  \frenchspacing\@vobeyspaces \@xverbatim}
\makeatother
\makeatletter
  \def\dgt@k{\dg@DX=-3 \dg@DY=2 \dg@SIZE=3} 
\makeatother

\makeatletter
  \def\dgt@kk{\dg@DX=3 \dg@DY=-1 \dg@SIZE=3}%
\makeatother

\theoremstyle{plain}
\newtheorem{thm}{Theorem}[section]
\newtheorem{cor}[thm]{Corollary}
\newtheorem{lem}[thm]{Lemma}
\newtheorem{pro}[thm]{Proposition}

\theoremstyle{definition}

\numberwithin{equation}{section}

\newcounter{rmnum}

\def\symbolnote#1#2{\let\thefootn=\thefootnote%
\renewcommand{\thefootnote}{\fnsymbol{footnote}}%
\footnotemark[#1]%
\footnotetext[#1]{#2}%
\let\thefootnote=\thefootn
}

\newfont{\bbb}{msbm10 scaled \magstep1}
\newfont{\bbc}{msbm8 scaled \magstep0}

\newcommand{\R}{\mbox{\bbb R}}

\newcommand{\N}{\mbox{\bbb N}}

\newcommand{\uin}{\mbox{\bbb I}}

\begin{document}


\title{On regularly branched maps
\newline
\newline
\newline
\it{dedicated to Professor S. Nedev for his 60th birthday}}

\author{H. Murat Tuncali}
\address{Department of Computer Science and Mathematics,
Nipissing University,
100 College Drive, P.O. Box 5002, North Bay, ON, P1B 8L7, Canada}
\email{muratt@nipissingu.ca}
\thanks{The first author was partially supported by his NSERC grant ?}

\author{Vesko Valov}
\address{Department of Computer Science and Mathematics, Nipissing University,
100 College Drive, P.O. Box 5002, North Bay, ON, P1B 8L7, Canada}
\email{veskov@nipissingu.ca}
\thanks{The second author was partially supported by his NSERC grant RGPIN 261914}

 

\begin{abstract}
Let $f\colon X\to Y$ be a perfect map between finite-dimensional metrizable spaces and $p\geq 1$.  It is shown that the space $C^*(X,\R^p)$ of all bounded maps from $X$ into $\R^p$ with the source limitation topology contains a dense $G_{\delta}$-subset consisting of $f$-regularly branched maps. Here, a map $g\colon X\to\R^p$ is $f$-regularly branched if, for every $n\geq 1$, the dimension of the set 
$\{z\in Y\times\R^p: |(f\times g)^{-1}(z)|\geq n\}$  is $\leq n\cdot\big(\dim f+\dim Y\big)-(n-1)\cdot\big(p+\dim Y\big)$. This is a parametric version of the Hurewicz theorem on regularly branched maps.
\end{abstract}

\maketitle

\markboth{H. M.~Tuncali and V.~Valov}{On regular branched maps}


\section{Introduction}

All spaces are assumed to be metrizable and all maps continuous. Moreover,
the function spaces in this paper, if not explicitely stated otherwise, are equipped with the source limitation topology.   
The paper is devoted to a parametric version of the Hurewicz theorem \cite{wh:33} on regularly branched maps. Recall that a map $g\colon X\to Z$ is called regularly branched (this term was introduced by Dranishnikov, Repov\v{s} and \v{S}\v{c}epin \cite{drs:91}) if $\dim B_n(g)\leq n\cdot\dim X-(n-1)\cdot\dim Z$ for any $n\geq 1$, where $B_n(g)=\{z\in Z: |g^{-1}(z)|\geq n\}$. 

\medskip\noindent
{\bf Hurewicz's Theorem}.  {\em Let $X$ be a finite-dimensional compactum and $p\geq 1$. Then the set of all regularly branched maps $g\colon X\to\R^p$ contains a dense $G_{\delta}$-subset of the space $C(X,\R^p)$.}  

\smallskip
We say that a map $g\colon X\to Z$ is regularly branched with respect to a fixed map  $f\colon X\to Y$ (briefly, $f$-regularly branched) if 

\medskip\noindent
$\dim B_n(f\times g)\leq n\cdot\big(\dim f+\dim Y\big)-(n-1)\cdot\big(\dim Z+\dim Y\big)$ for every $n\geq 1$,

\medskip\noindent
where $\dim f=\sup\{\dim f^{-1}(y): y\in Y\}$.  Obviously, when $f$ is a constant map, i.e.,  $Y$ is a point, the notions of $f$-regularly branched and regularly branched maps coincide.  Next theorem is our main result. 

\begin{thm}
Let $f\colon X\to Y$ be a $\sigma$-perfect map between finite-dimensional spaces and $p\geq 1$. Then the space $C^*(X,\R^p)$ contains a dense $G_{\delta}$-subset $\mathcal{H}$ consisting of $f$-regularly branched maps.
\end{thm}

Here, $C^*(X,\R^p)$ is the set of all bounded maps from $X$ into $\R^p$ and $f$ is said to be $\sigma$-perfect if $X$ is the union of  its closed subsets $X_i$, $i=1,2,..$,  such that $f(X_i)\subset Y$ are closed and each restriction $f|X_i$ is perfect.  

\begin{cor}
Let the integers $k$, $p$, $m$ and $n$ satisfy the inequality $k+m+1\leq (p-k)n$. Then, for any $\sigma$-perfect map $f\colon X\to Y$ with $\dim f\leq k$ and $\dim Y\leq m$, the space 
$C^*(X,\R^p)$ contains a dense $G_{\delta}$-subset of maps $g$ such that 
$ |(f\times g)^{-1}(z)|\leq n$ for every $z\in Y\times\R^p$.
\end{cor}   

Corollary 1.2 follows directly from Theorem 1.1. Indeed,  under the hypotheses of this corollary, if $g\in C^*(X,\R^p)$ is $f$-regularly branched, then $\dim B_{n+1}(f\times g)\leq (n+1)(k+m)-n(p+m)\leq -1$. So, $f\times g$ is $\leq n$-to-one for all $f$-regularly branched maps. Let us also mention next corollary of Theorem 1.1 (it follows, actually, from Corollary 1.2) established by the authors in \cite{tv1} and providing positive solutions of two hypotheses of Bogatyi-Fedorchuk-van Mill \cite{bfm:00}.
 
\begin{cor}
Let $f\colon X\to Y$ be a $\sigma$-perfect map with $\dim f\leq k$ and $\dim Y\leq m$.  Then, for any $p\geq 1$, $C^*(X,\R^{p+k})$ contains a dense $G_{\delta}$-subset consisting of maps $g$ such that $|(f\times g)^{-1}(z)|\leq \max\{k+m-p+2,1\}$ for all $z\in Y\times\R^p$.
\end{cor}   

If  $p\geq 2k+m+1$, then Corollary 1.2 (as well as, Corollary 1.3) yields  the existence of a dense and $G_{\delta}$-subset $G$ of $C^*(X,\R^p)$ such that $f\times g$ is one-to-one for every $g\in G$.  Hence, all $f\times g$, $g\in G$, are embeddings provided $f$ is a perfect map.   So, we obtain a    
parametric version of the N\"{o}beling-Pontryagin embedding theorem which was established in \cite{bp:96}, \cite{bp:98} and \cite{tv4}. 

The question if the set $\mathcal{H}$ from Theorem 1.1 can consist of maps $g$ such that  
$\dim B_n(f\times g)\leq n\cdot\dim X-(n-1)\cdot\big(p+\dim Y\big)$ for every $n\geq 1$
was raised in the first version of this paper. The reviewer and S. Bogatyi independently provided a negative answer. Here is the example suggested by Bogatyi:\\
Let $T$ be a metrizable compactum not embeddable in $\R^{2m}$, $m\geq 2$, such that $\dim T\leq m$. Take the disjoint sum $X=\uin^m\oplus T$ and the map $f\colon X\to\uin^m$, 
$f(x)=x$ if $x\in\uin^m$ and $f(x)=x_0\in\uin^m$ if $x\in T$. The existence of a map $g\colon X\to\R^{m+2}$ with the above property would imply  
that $g$ embeds $T$ into $\R^{m+2}$ which is impossible because $m+2\leq 2m$.

\smallskip
Let us also note that, by \cite[Corollary 11]{b}, for every $m$ there exists a polyhedron $X$ with $\dim X=m$ such that every map $g\in C(X,\R^{m+1})$ has a fiber containing at least $m+1$ points. Therefore, the inequality in the definition of a regularly branched map
$\dim B_n(f\times g)\leq n\cdot\big(\dim f+\dim Y\big)-(n-1)\cdot\big(\dim Z+\dim Y\big)$  can not be improved.   

\medskip\noindent
The original proof of Theorem 1.1 was quite complicated. Based on our previous results from \cite{tv:02} and \cite{tv1}, the referee of this paper found very elegant proof of Theorem 1.1 and this proof is presented here. Moreover, we 
provide a unified method for proving the results used in the proof of Theorem 1.1. 
This method is extracted from our previous papers \cite{tv:02}, \cite{tv:03}, \cite{tv1} and
\cite{tv4}. Is is based on selection theorems established by the second author and V. Gutev in
\cite{gv:99} and \cite{gv:02}.

\smallskip
We acknowledge the referee's proof of Theorem 1.1. We are also grateful to S. Bogatyi who provided the described above counterexample to our question.  


\section{Some preliminary results}

\bigskip
First, we provide some information about the source limitation topology.  This topology can be described as follows: If $(M,d)$ is a metric space, then 
a set $U\subset C(X,M)$ is open if for every $g\in U$ there exists 
a continuous function $\alpha\colon X\to (0,\infty)$ such that $\overline{B}(g,\alpha)\subset U$. Here, $\overline{B}(g,\alpha)$ denotes the set  
$\{h\in C(X,M):d(g(x),h(x))\leq\alpha (x)\hbox{}~~\mbox{for each 
$x\in X$}\}$. 
The source limitation topology doesn't depend on the metric $d$ if $X$ is paracompact \cite{nk:69}. Moreover, $C(X,M)$ with this topology has the Baire property provided $(M,d)$ is a complete metric space \cite{jm:75}. Obviously, the source limitation topology coincides with the uniform convergence topology generated by $d$ in case $X$ is compact.  One can show that $C^*(X,\R^p)$ is open in $C(X,\R^p)$ with respect to the source limitation topology when the Euclidean metric on $\R^p$ is considered.  Therefore, $C^*(X,\R^p)$ equipped with this topology also has the Baire property.

\bigskip
We are going to establish a background of the general method discussed in the introduction.  
Throughout this section $K$ is a closed and convex subset  of a given Banach space $E$ and   $f\colon X\to Y$  a  perfect surjective map between metrizable spaces.
Suppose, for every $y\in Y$, we are given a set ${\mathcal C}(y)\subset C^*(X,K)$ such that if $h\in C^*(X,K)$ and $h|f^{-1}(y)=g|f^{-1}(y)$ for some $g\in\mathcal{C}(y)$, then
$h\in\mathcal{C}(y)$. The last property means that the set $\mathcal{C}(y)$ is determined by the restrictions $g|f^{-1}(y)$. That is why, sometimes, we consider $\mathcal{C}(y)$ as a class of functions on $f^{-1}(y)$. Let 
$\displaystyle {\mathcal C}(H)=\bigcap_{y\in H}{\mathcal C}(y)$, where $H\subset Y$.  
We also consider the set-valued map $\displaystyle\psi\colon Y\to C^*(X,K)$, defined by  
$\psi(y)=C^*(X,K)\backslash {\mathcal C}(y)$.

\begin{lem}
Suppose for every $y\in Y$ and every $g\in {\mathcal C}(y)$ there exists a neighborhood $V_y$ of $y$ in $Y$ and $\delta_y>0$ such that $h\in\mathcal{C}(V_y)$ provided $h|f^{-1}(V_y)$ is $\delta$-close to $g|f^{-1}(V_y)$. Then $\mathcal{C}(Y)$ is open in $C^*(X,K)$.   
Moreover, $\psi$ has a closed graph when $C^*(X,K)$ is equipped with the uniform convergence topology.
\end{lem}

\begin{proof}
We follow some ideas from \cite{bv}.
Let $(y_0,g_0)\in Y\times C^*(X,K)\backslash G_{\psi}$, where $C^*(X,K)$ possesses the uniform convergence topology and $G_{\psi}$ is the graph of $\psi$. Hence, $g_0\not\in\psi(y_0)$, so $g_0\in\mathcal{C}(y_0)$. Take $V_{y_0}$ and $\delta_{y_0}>0$ satisfying the hypotheses of the lemma, and let $W$ denote the $\delta_{y_0}$-neighborhood of $g_0$ in $C^*(X,K)$. Then $V_{y_0}\times W$ is a neighborhood of $(y_0,g_0)$ disjoint from $G_{\psi}$. Thus, $G_{\psi}$ is closed.

To show that $\mathcal{C}(Y)$ is open in $C^*(X,K)$ with respect to the source limitation topology, we fix $g_0\in\mathcal{C}(Y)$. Since, for every $y\in Y$, $g_0\in\mathcal{C}(y)$, we choose neighborhoods $V_y$ and positive numbers $\delta_y\leq 1$ satisfying the conditions of the lemma. We can assume that $\{V_y:y\in Y\}$ is a locally finite cover of $Y$, and consider the set-valued map $\varphi\colon Y\to (0,1]$, $\varphi(y)=\cup\{(0,\delta_z]:y\in V_z\}$. 
Then, by \cite[Theorem 6.2]{rs}, $\varphi$ admits a continuous selection $\beta\colon Y\to (0,1]$, and let $\alpha=\beta\circ f$. It remains only to show that if $g\in C^*(X,K)$ with $d\big(g_0(x),g(x)\big)<\alpha(x)$ for all $x\in X$, where $d$ is the metric on $E$ generated by its norm, then $g\in\mathcal{C}(Y)$. So, we take such a $g$ and fix $y\in Y$. 
Then, there exists $z\in Y$ with $y\in V_{z}$ and such that $\alpha(x)\leq\delta_{z}$ for all $x\in f^{-1}(y)$. Now, select a map $h\in C^*(X,K)$ coinciding with $g$ on $f^{-1}(y)$ and satisfying the inequality $d\big(h(x),g_0(x)\big)\leq\delta_z$ for each $x\in X$. According to the choice of $V_z$, $h\in\mathcal{C}(y)$. Hence, $g\in\mathcal{C}(y)$ because
$g|f^{-1}(y)=h|f^{-1}(y)$. Therefore, $\mathcal{C}(Y)$ is open in $C^*(X,K)$.       
\end{proof}

Recall that a closed subset $F$ of the metrizable space $M$ is said to be a $Z_m$-set in $M$,   if the set $C(\uin^m,M\backslash F)$ is dense in $C(\uin^m,M)$ with respect to the uniform convergence topology, where $\uin^m$ is the $m$-dimensional cube. If  
$F$ is a $Z_m$-set in $M$ for every $m\in\N$, we say that $F$ is a $Z$-set in $M$.

\begin{lem}
Let $y\in Y$ and ${\mathcal C}(y)$, considered as a subset of $C(f^{-1}(y),K)$, satisfy the following condition:
\begin{itemize}
\item For every $k\in\N$ $($resp., $k=m$$)$ the set of all maps $h\in C(\uin^k\times f^{-1}(y),K)$ with $h|(\{z\}\times f^{-1}(y))\in\mathcal{C}(y)$ for each $z\in\uin^k$, is dense in $C(\uin^k\times f^{-1}(y),K)$ with respect to the uniform convergence topology.
\end{itemize}
Then, for every $\alpha\colon X\to (0,\infty)$ and $g\in C^*(X,K)$, $\psi(y)\cap\overline{B}(g,\alpha)$  is a $Z$-set $($resp., $Z_m$-set$)$ in 
$\overline{B}(g,\alpha)$ provided $\overline{B}(g,\alpha)$ is considered as a subset of $C^*(X,K)$ equipped with the uniform convergence topology and
$\psi(y)\subset C^*(X,K)$ is closed.
\end{lem}

\begin{proof}
See the proof of \cite[Lemma 2.8]{tv:02}.
\end{proof}

\begin{lem}
Let  $Y$ be a $C$-space $($resp., $\dim Y\leq m$$)$ and  the family $\{{\mathcal C}(y)\}_{ y\in Y}$ satisfies the following conditions:

\begin{itemize}
\item[(a)] the map $\psi$ has a closed graph;
\item[(b)] $\psi(y)\cap\overline{B}(g,\alpha)$ is a $Z$-set $($resp., $Z_m$-set$)$ in $\overline{B}(g,\alpha)$
for any continuous function $\alpha\colon X\to (0,\infty)$, $y\in Y$ and $g\in C^*(X,K)$,  where 
$\overline{B}(g,\alpha)$ is considered as a subspace of $C^*(X,K)$ with the uniform convergence topology. 
\end{itemize}
Then $\mathcal{C}(Y)$ is dense in $C^*(X,K)$ with respect to the source limitation topology.
\end{lem}

\begin{proof}
It suffices to show that, for fixed $g_0\in C^*(X,K)$ and a continuous function $\alpha\colon X\to (0,\infty)$, there exists $g\in \overline{B}(g_0,\alpha)\cap\mathcal{C}(Y)$. We equip $C^*(X,K)$ with the uniform convergence topology and 
consider
the constant convex-valued map $\phi\colon Y\to C^*(X,K)$, 
$\phi(y)=\overline{B}(g_0,\alpha_1)$, where $\alpha_1(x)=\min\{\alpha (x), 1\}$.
Because of the conditions (a) and (b), we  
can apply the selection theorem \cite[Theorem 1.1]{gv:99} (resp., \cite[Theorem 1.1]{gv:02}) to obtain a  continuous map $h\colon Y\to C^*(X,K)$ such that  $h(y)\in\phi(y)\backslash\psi(y)$ for every $y\in Y$.
Observe that $h$ is a map from $Y$ into $\overline{B}(g_0,\alpha_1)$ such that $h(y)\in\mathcal{C}(y)$ for every $y\in Y$. Then
$g(x)=h(f(x))(x)$, $x\in X$, defines a bounded map $g\in \overline{B}(g_0,\alpha)$ such that  $g|f^{-1}(y)=h(y)|f^{-1}(y)$, $y\in Y$. Therefore, $g\in\mathcal{C}(y)$ for all $y\in Y$, i.e., $g\in \overline{B}(g_0,\alpha)\cap\mathcal{C}(Y)$.
\end{proof}


\section{Finite-to-one maps}

\bigskip 
In this section 
we provide a non-compact version, see Proposition 3.1 below, of the Levin-Lewis result \cite[Proposition 4.4]{ll:02}. Note that, for separable metrizable spaces, Proposition 3.1 follows from \cite[Lemma 2]{ht:85}. 

\begin{pro}
Let $f\colon X\to Y$ be a perfect $0$-dimensional map with $\dim Y\leq m$. Then $C^*(X)$
contains a dense $G_{\delta}$-subset of  maps $g$ with  each fiber of $f\times g$ containing at most $m+1$ points.
\end{pro}

\begin{proof}
We take a map $\theta\colon X\to Q$ such that $f\times\theta\colon X\to Y\times Q$ is an embedding
 (such a $\theta$ exists by \cite{bp:98} or \cite{tv4}) with $Q$ being the Hilbert cube, a countable base $\{W_i\}_{i\in\N}$ of open sets in $Q$. Let $\mathcal A$ be the collection of the closures of $\theta^{-1}(W_i)$ in $X$, $i\geq 1$.  There are countably many families $\Gamma=\{A_1,A_2,..,A_{m+2}\}$ consisting of $m+2$ disjoint 
elements of $\mathcal A$. 
For any such $\Gamma$ and $y\in Y$ let  $\mathcal{C}_{\Gamma}(y)$ denote the set of all $g\in C^*(X)$ such that each $g^{-1}(z)\cap (f^{-1}(y)$, $z\in\R$, meets at most $m+1$ elements of $\Gamma$. Following Section 2, for $H\subset Y$, let $\mathcal{C}_{\Gamma}(H)=\cap\{\mathcal{C}_{\Gamma}(y):y\in H\}$. Since 
the intersection of all $\mathcal{C}_{\Gamma}(Y)$ consists of maps $g$ such that each fiber of $f\times g$ contains at most $m+1$ points, it suffices to show that any 
$\mathcal{C}_{\Gamma}(Y)$ is open and dense in $C^*(X)$.

\begin{lem}
Let $\Gamma=\{G_1,..,G_{m+2}\}$ and $y\in Y$ be fixed. Then, for every $g\in {\mathcal C}_{\Gamma}(y)$ there exists a neighborhood $V$ of $y$ in $Y$ and $\delta>0$ such that $h\in\mathcal{C}_{\Gamma}(V)$ provided $h|f^{-1}(V)$ is $\delta$-close to $g|f^{-1}(V)$.
\end{lem} 

\begin{proof}
Assume this is not true for some $g_0\in\mathcal{C}_{\Gamma}(y)$. Then, there exist neighborhoods $V_i$, $i\geq 1$, of $y$ in $Y$, functions $g_i\in C^*(X)$,  points $y_i\in V_i$ and $z_i\in\R$ such that $g_i|f^{-1}(V_i)$ is $\displaystyle 1/i$-close to $g_0|f^{-1}(V_i)$ but $g_i^{-1}(z_i)\cap f^{-1}(y_i)$ meets all $\leq m+2$ elements of $\Gamma$. 
Since $f$ is closed, we can suppose that $U_i=f^{-1}(V_i)\subset g_0^{-1}(W_i)$ with $U_i$ and $W_i$ being $\displaystyle 1/i$ neighborhoods of $f^{-1}(y)$ and $g_0(f^{-1}(y))$ in $X$ and $\R$, respectively, and $z_i\in W_i$. Passing to subsequences, we may also suppose that $\lim z_i=z_0\in g_0(f^{-1}(y))$. Then 
$g_0^{-1}(z_0)\cap f^{-1}(y)$ intersects at most $m+1$ elements of $\Gamma$, let say the first $m+1$. 
Take points $a_i\in g_i^{-1}(z_i)\cap f^{-1}(y_i)$ and $b_i\in f^{-1}(y)$ such that $a_i\in G_{m+2}$ and $dist(a_i,b_i)\leq 1/i$ for all $i$. Again, we can assume that $\lim b_i=b_0$ for some $b_0\in f^{-1}(y)$. Then $\lim a_i=b_0\in g_0^{-1}(z_0)\cap f^{-1}(y)$, so $b_0\not\in G_{m+2}$. This implies that $a_i\not\in G_{m+2}$ for almost all $i$ which contradicts the choice of the points $a_i$.  
\end{proof} 
 
Therefore, combining Lemma 3.2 and Lemma 2.1, we may conclude that each $\mathcal{C}_{\Gamma}(Y)$ is open in $C^*(X)$ and the set-valued map $\psi_{\Gamma}\colon Y\to C^*(X)$, $\psi_{\Gamma}(y)=C^*(X)\backslash\mathcal{C}_{\Gamma}(y)$, has a closed graph when $C^*(X)$ carries the uniform convergence topology.

\begin{lem}
For any $\Gamma$ and $y\in Y$, the set of all functions 
$g\in C\big(\uin^m\times f^{-1}(y)\big)$ such that $g|\big(\{z\}\times f^{-1}(y)\big)\in\mathcal{C}_{\Gamma}(y)$ for each $z\in\uin^m$, is dense in  $C\big(\uin^m\times f^{-1}(y)\big)$ 
\end{lem}

\begin{proof}
By the Levin-Lewis result \cite[Proposition 4.4]{ll:02}, every $h\in C\big(\uin^m\times f^{-1}(y)\big)$ cam be approximated by functions $g\in C\big(\uin^m\times f^{-1}(y)\big)$ such that each $g^{-1}(t)\cap \big(\{z\}\times f^{-1}(y)\big)$, $z\in\uin^m$ and $t\in\R$, contains at most $m+1$ points. This implies that $g|\big(\{z\}\times f^{-1}(y)\big)\in\mathcal{C}_{\Gamma}(y)$ for each $z\in\uin^m$, and we are done.
\end{proof}

Finally, the combination of Lemma 3.3 and Lemma 2.1 - 2.3, yields that every 
$\mathcal{C}_{\Gamma}(Y)$ is dense in $C^*(X)$. This completes the proof of Proposition 3.1.
\end{proof}



\section{Proof of Theorem 1.1}

\bigskip
One of the components of the proof of Theorem 1.1 is Theorem 4.1 below. It is a parametric version of the Hurewicz result \cite{wh:33} that every $n$-dimensional compactum admits a $0$-dimensional map into $\uin^n$. For finite-dimensional compact spaces this version was proved by Pasynkov \cite{bp:96} (announced in 1975). Torunczyk \cite{ht:85} also established such a theorem for finite-dimensional separable spaces. In the present form, Theorem 4.1 was obtained by the authors \cite{tv:02}. The proof presented here  follows the general method from Sections 2 and 3. Pasynkov's theorem, mentioned above, is also used, but we provide an easy proof of that theorem.

\begin{thm}
Let $f\colon X\to Y$ be a $\sigma$-perfect $n$-dimensional map with $Y$ being a $C$-space. Then, for every $0\leq k\leq n$, $C^*(X,\R^k)$ contains a dense $G_{\delta}$-subset of maps $g$ such that $f\times g$ is $(n-k)$-dimensional. 
\end{thm}   
 
\begin{proof}
It is easily seen that the proof is reduced to the case when $f$ is perfect. Following the general shem from Section 2,
for every $\epsilon>0$ and $y\in Y$, let $\mathcal{C}_{\epsilon}(y)$ be the set of all maps $g\in C^*(X,\R^k)$ satisfying the following condition: every set $f^{-1}(y)\cap g^{-1}(z)$, $z\in\R^k$, can be covered by a finite family $\gamma$ of open sets in $X$ each of diameter $\leq\epsilon$ and any point of $X$ is contained in at most $n-k+1$ elements of $\gamma$. We need to show that every $\mathcal{C}_{\epsilon}(Y)$ is open and dense in $C^*(X,\R^k)$.
The proof of next lemma is similar to that one of Lemma 3.2.

\begin{lem}
Let $\epsilon>0$ and $y\in Y$ be fixed. Then, for every $g\in {\mathcal C}_{\epsilon}(y)$ there exists a neighborhood $V$ of $y$ in $Y$ and $\delta>0$ such that $h\in\mathcal{C}_{\epsilon}(V)$ provided $h|f^{-1}(V)$ is $\delta$-close to $g|f^{-1}(V)$.
\end{lem} 

As above, Lemma 4.2 implies that  all ${\mathcal C}_{\epsilon}(Y)$ are open in $C^*(X,\R^k)$ and the set-valued map $\psi_{\epsilon}\colon Y\to C^*(X,\R^k)$, $\psi_{\epsilon}(y)=
C^*(X,\R^k)\backslash {\mathcal C}_{\epsilon}(y)$, has a closed graph when $C^*(X,\R^k)$ is equipped with the uniform convergence topology. 

The density of the sets ${\mathcal C}_{\epsilon}(Y)$ in $C^*(X,\R^k)$ follows from the lemma below:

\begin{lem}
For any $\epsilon>0$, $m\geq 1$ and $y\in Y$, the set of all maps 
$g\in C\big(\uin^m\times f^{-1}(y),\R^k\big)$ such that $g|\big(\{z\}\times f^{-1}(y)\big)\in\mathcal{C}_{\epsilon}(y)$ for each $z\in\uin^m$, is dense in  $C\big(\uin^m\times f^{-1}(y),\R^k\big)$ 
\end{lem}

\begin{proof}
The proof of this lemma is the same as the proof of Lemma 3.3. The only difference now is that, instead of the Levin-Lewis theorem, we use the Pasynkov result formulated in Proposition 4.4 below. 
\end{proof}

Combining all lemmas in Section 2 and Section 3, we can complete the proof of Theorem 4.1.
Therefore, we need only to provide a proof of Proposition 4.4.    

\begin{pro}
Let $K$ be a compactum of dimension $\leq n$ and $0\leq k\leq n$. Then, for every $m\geq 1$, the set of all maps $g\in C\big(\uin^m\times K,\R^k\big)$ such that   
$\pi\times g$ is   
$(n-k)$-dimensional, is dense in $C\big(\uin^m\times K,\R^k\big)$ $($here $\pi\colon\uin^m\times K\to\uin^m$ denotes the projection$)$.
\end{pro}

Observe that the validity of the case $k=n$ implies the validity of all other cases. Indeed,
if $h\in C(\uin^m\times K,\R^k)$ and $\eta>0$, we lift $h$ to a map $h_1\colon\uin^m\times K\to\R^n$ such that $h=p\circ h_1$, where $p\colon\R^n\to\R^k$ is the canonical projection. Next, take $g_1\in C(\uin^m\times K,\R^n)$ $\eta$-close to $h_1$ and such that $\pi\times g_1$ is $0$-dimensional. Then,  $g=p\circ g_1$ is $\eta$-close to $h$ and $\pi\times g$ is   
$(n-k)$-dimensional. So, we can suppose that $k=n$.

Since $\dim K\leq n$, by the Hurewicz theorem \cite{wh:33}, there exists a $0$-dimensional map $g\colon K\to\uin^n$. Then $\pi\times\overline{g}$, where $\overline{g}$ is the composition of the projection $\pi_K\colon\uin^m\times K\to K$ and $g$, is also $0$-dimensional. According to \cite[$(ii)\Leftrightarrow (iii)$]{l}, almost all maps $g\in C(\uin^m\times K,\R^n)$ have the property $\dim (\pi\times g)\leq 0$. This completes the proof of Proposition 4.4. Finally, let us note that Levin's result 
\cite[$(ii)\Leftrightarrow (iii)$]{l}, which was used in this proof, has a very short proof. 
As a result, we obtain a proof of Proposition 4.4 which is quite easier than the original one from
\cite{bp:96}.
\end{proof}

\bigskip\noindent
{\em Proof of Theorem $1.1$} 
Let show first that the proof of Theorem 1.1 can be reduced to the case $f$ is perfect. 
Suppose $X$ is the union of an increasing sequence of its closed sets $X_i$ such that each restriction $f_i=f|X_i$ is perfect with $Y_i=f(X_i)\subset Y$ being closed.  Then, applying Theorem 1.1 for every map $f_i\colon X_i\to Y_i$, and using that 
the maps $\pi_i\colon C^*(X,\R^p)\to C^*(X_i,\R^p)$, $\pi_i(g)=g|X_i$, are surjective and open, we conclude that 
there exists a dense $G_{\delta}$-set $G\subset C^*(X,\R^p)$ consisting of maps $g$ such that $g_i=g|X_i$ is $f_i$-regularly branched for every $i$.  Let $g\in G$ and $n\geq 1$.  
For any $i$ the  set $B_n(f_i\times g_i)$ is $F_{\sigma}$ in $(f_i\times g_i)(X_i)$ \cite{re:95} and  $(f_i\times g_i)(X_i)\subset Y\times\R^p$ is closed (recall that each $Y_i\subset Y$ is closed and 
the map $f_i\times g_i\colon X_i\to Y_i\times\R^p$ is perfect).  
So, all of the sets $B_n(f_i\times g_i)$ are $F_{\sigma}$ in $Y\times\R^p$.  Moreover, 
 $\dim B_n(f_i\times g_i)\leq n\cdot\big(\dim f_i+\dim Y_i\big)-(n-1)\cdot\big(p+\dim Y_i\big)\leq n\cdot\big(\dim f+\dim Y\big)-(n-1)\cdot\big(p+\dim Y\big)$.
Therefore, $\dim \bigcup_{i=1}^{\infty}B_n(f_i\times g_i)\leq n\cdot\big(\dim f+\dim Y\big)-(n-1)\cdot\big(p+\dim Y\big)$. On the other hand, 
$B_n(f\times g)\subset\bigcup_{i=1}^{\infty}B_n(f_i\times g_i)$.  Consequently,
$\dim B_n(f\times g)\leq n\cdot\big(\dim f+\dim Y\big)-(n-1)\cdot\big(p+\dim Y\big)$ for every $g\in G$ and $n\geq 1$. Hence, $G$ consists of $f$-regularly branched maps. Thus, everywhere below we may assume that $f$ is perfect. Moreover, we can also assume that $p>\dim f$ because, according to the definition, every  $g\in C(X,\R^p)$ is $f$-regularly branched provided $p\leq\dim f$.

The remaining part of the proof, presented below, was suggested by the referee of this paper.    
 
It is easily seen that, by Theorem 4.1, we can assume $\dim f=0$. So, everywhere below $f$ is a perfect $0$-dimensional map, $p\geq 1$ and $\dim Y=m$. Let $l=l(m,p)=[m/p]+1$, where $[m/p]$ denotes the integer part of $m/p$. 

We show by induction on $p$ that $f\times g$ is at most $l$-to-$1$ for almost all maps $g\in C^*(X,\R^p)$. For $p=1$, it follows from Proposition 3.1. Assume $p>1$ and let $m=(l-1)p+t$, $0\leq t<p$. Decompose $Y=Y_1\cup Y_2$ such that $Y_1$ is an $F_{\sigma}$-subset of $Y$ with $\dim Y_1\leq m-l=(l-1)(p-1)+t-1$ and $\dim Y_2\leq l-1$. Let also $g=g_1\times g_2\colon X\to \R^{p-1}\times\R$. 
Since $\displaystyle [(m-l)/(p-1)]+1=l$, according to the induction hypothesis, $g_1$ can be approximated by a map $g_1^*\colon X\to\R^{p-1}$ such that $f\times g_1^*$ is at most $l$-to-$1$ on  $f^{-1}(Y_1)$. Denote by $B$ the union of all fibers of $f\times g_1^*$ having more than $l$ points. Then $B$ is $F_{\sigma}$ in $X$ and disjoint from $f^{-1}(Y_1)$, so 
$f(B)\subset Y_2$. Once again by induction hypothesis, $g_2$ can be approximated by a map $g_2^*\colon X\to\R$ such that $f\times g_2^*$ is at most $l$-to-$1$ on $f^{-1}(f(B))$. Thus, $g$ can be approximated by the map $g^*=g_1^*\times g_2^*$ such that $f\times g^*$ is at most $l$-to-$1$. This implies that the maps $g\in C^*(X,\R^p)$ such that $f\times g$ is at most $l$-to-$1$ form a dense subset of $C^*(X,\R^p)$. To complete the induction, we need to show that this set is also $G_{\delta}$ in $C^*(X,\R^p)$. To this end, following the proof of Proposition 3.1, we  take a map $\theta\colon X\to Q$ such that $f\times\theta\colon X\to Y\times Q$ is an embedding, and a
countable base $\{W_i\}_{i\in\N}$ of open sets in $Q$. We also consider the collection $\mathcal A$ of all closures of $\theta^{-1}(W_i)$ in $X$, $i\geq 1$.  There are countably many families $\Gamma=\{A_1,A_2,..,A_{l+1}\}$ consisting of $l+1$ disjoint 
elements of $\mathcal A$ and 
for any such $\Gamma$ and $y\in Y$ let  $\mathcal{C}_{\Gamma}(y)$ denote the set of all $g\in C^*(X,\R^p)$ such that each $g^{-1}(z)\cap f^{-1}(y)$, $z\in\R^p$, meets at most $l$ elements of $\Gamma$. As in Section 3, we can show that any set 
$\mathcal{C}_{\Gamma}(Y)=\cap\{\mathcal{C}_{\Gamma}(y):y\in Y\}$ is open in $C^*(X,\R^p)$.
Therefore, the maps $g\in C^*(X,\R^p)$ with $f\times g$ being at most $l$-to-$1$ form a $G_{\delta}$-set in $C^*(X,\R^p)$ as the intersection of all $\mathcal{C}_{\Gamma}(Y)$.

Now, we can finish the proof of Theorem 1.1. Let $Y_i\subset Y$, $0\leq i\leq m$, be $F_{\sigma}$-subsets of $Y$ such that $Y_0\subset Y_1\subset...\subset Y_m$, $\dim Y_i\leq i$ and $\dim Y\backslash Y_i\leq m-i-1$. Then, from what we proved above, it follows that $C^*(X,\R^p)$ contains a dense $G_{\delta}$-subset $G$ of maps $g$ such that $f\times g$ is at most $l(i,p)$-to-$1$ on $f^{-1}(Y_i)$ for every $0\leq i\leq m$. Moreover, in addition, we may require by \cite{tv:03} that $g(f^{-1}(y))$ is $0$-dimensional for all $y\in Y$ and all $g\in G$. It remains only to show that every $g\in G$ is $f$-regularly branched. So, we fix $g\in G$ and $n\geq 1$, and let $\pi_Y\colon Y\times\R^p\to Y$ be the projection onto $Y$. Since 
$B_n(f\times g)$ is $F_{\sigma}$ in $(f\times g)(X)$ and $\pi_Y|(f\times g)(X)$ is a perfect map, $\pi_Y\big(B_n(f\times g)\big)$ is $F_{\sigma}$ in $Y$. Moreover, since 
each $g(f^{-1}(y))$ is $0$-dimensional, $\dim B_n(f\times g)$ is at most the dimension of $\pi_Y\big(B_n(f\times g)\big)$. On the other hand, if $(f\times g)^{-1}(y,z)$ contains $\geq n$ points, then $y\not\in Y_{p(n-1)-1}$. Hence, 
$\pi_Y\big(B_n(f\times g)\big)$ is contained in $Y\backslash Y_{p(n-1)-1}$. Consequently,
$\dim \pi_Y\big(B_n(f\times g)\big)\leq m-(n-1)p$, so is $\dim B_n(f\times g)$. Since
 $n(\dim f+\dim Y)-(n-1)(p+\dim Y)=m-(n-1)p$, the last inequality shows that $g$ is regularly $f$-branched.

\bigskip


\begin{thebibliography}{99}

\bibitem{b}
S.~Bogatyi, {\em The geometry of maps into Euclidean space}, Uspekhi Mat. Nauk {\bf 53:5} (1998), 27--56 (in Russian); English translation: Russian Math. Surveys {\bf 53:5} (1998), 893--920.

\bibitem{bfm:00}
S.~Bogatyi, V.~Fedorchuk and J. van Mill, {\em On mappings of compact spaces into Cartesian spaces}, Topology and Appl. {\bf 107} (2000), 13--24.

\bibitem{bv}
S.~Bogatyi and V.~Valov, {\em Roberts' type theorems on dimension}, submitted.
 
\bibitem{ch:96}
A.~Chigogidze, {\em Inverse Spectra}, North Holland, Amsterdam, 1996.

\bibitem{drs:91}
A.~Dranishnikov, D. ~Repov\v{s} and E.~\v{S}\v{c}epin, {\em On intersections of compacta of complementary dimensions in Euclidean space}, Topology Appl. 38 (1991), 237--253.

\bibitem{re:95}
R.~Engelking, {\em Theory of dimensions: Finite and Infinite}
(Heldermann Verlag, Lemgo, 1995).

\bibitem{gv:99}
V.~Gutev and V.~Valov, {\em Dense families of selections and finite-dimensional spaces}, Set-Valued Analysis {\bf 11} (2003), 373--391.

\bibitem{gv:02}
V.~Gutev and V.~Valov, {\em Continuous selections and $C$-spaces}, Proc. Amer.
Math. Soc. {\bf 130, 1} (2002), 233--242.
 
\bibitem{wh:33}
W.~Hurewicz, {\em \"{U}ber Abbildungen von endlichdimensionalen R\"{a}umen auf Teilmengen Cartesischer 
R\"{a}ume}, Sgb. Preuss. Akad., {\bf 34} (1933), 754--768.

\bibitem{nk:69}
N.~Krikorian, \emph{A note concerning the fine topology on function spaces},
  Compos. Math. \textbf{21} (1969), 343--348.

\bibitem{ll:02}
M.~Levin and W. Lewis, {\em Some mapping theorems for extensional dimension}, Israel J. Math. {\bf 133} (2003), 61--76.

\bibitem{l}
M.~Levin, {\em Bing maps and finite-dimensional maps}, Fund. Math. {\bf 151} (1996), 47--52.

\bibitem{jm:75}
J.~Munkers, {\em Topology} (Prentice Hall, Englewood Cliffs, NY, 1975).

\bibitem{bp:96}
B.~Pasynkov, {\em On geometry of continuous maps of finite-dimensional compact metric spaces}, Trudy Mat. Inst. Steklova
{\bf 212} (1996), 147--172 (in Russian); English translation: Proc. Steklov Inst. Math. {\bf 212:1} (1996), 138--162.

\bibitem{bp:98}
\bysame, {\em On geometry of continuous maps of countable functional weight}, Fundam.
Prikl. Matematika {\bf 4, 1} (1998), 155--164 (in Russian).

\bibitem{rs}
D.~Repov\v{s} and P.~Semenov, {\em Continuous selections of multivalued mappings} (Math. and its Appl. {\bf 455}, Kluwer, Dordrecht, 1998).

\bibitem{ht:80}
H.~Torunczyk, {\em On $CE$-images of the Hilbert cube and characterization of $Q$-manifolds}, Fund. Math. {\bf 106} (1980), 31--40.

\bibitem{ht:85}
\bysame, {\em Finite-to-one restrictions of continuous functions}, Fund. Math.
{\bf 125} (1985), 237--249.

\bibitem{tv:02}
M.~Tuncali and V.~Valov, {\em On dimensionally restricted maps}, Fund. Math. 175, 1 (2002), 35--52.

\bibitem{tv:03}
\bysame, \bysame, {\em On finite-dimensional maps II}, Topology and Appl. {\bf 132} (2003), 81--87.
 
\bibitem{tv1}
\bysame, \bysame, {\em On finite--to-one-maps}, Canadian Math. Bulletin, to appear.

\bibitem{tv4}
\bysame, \bysame, {\em On finite-dimensional maps}, Tsukuba J. Math. {\bf 28, 1} (2004), 155--167.

\end{thebibliography}
\end{document}